\documentclass[psamsfonts,11pt]{article}
\usepackage{amsfonts,amssymb,eucal}
\usepackage{amsmath,amsthm}
\catcode `\@=11 \@addtoreset{equation}{section}

\newtheorem{cor}{Corollary}[section]
\newtheorem{lem}{Lemma}

\newtheorem{thm}{Theorem}

\input xy
\xyoption{all}

\title{A compact group action which raises dimension to infinity}

\author{Zhiqing Yang\\ email:\ zyang@math.ucla.edu}

\begin{document}
\maketitle{}

\section{Introduction}

\bigskip
Interest in group actions with the property that the orbit space
has greater dimension than the original space is motivated in part
of P.A.Smith's generalization of the Hilbert's fifth problem. It
asserts that among all locally compact groups only Lie group can
act effectively on manifolds. It is now known as the Hilbert-Smith
conjecture. The Hilbert-Smith conjecture is known to be equivalent
to the conjecture that a compact topological group acting on a
compact topological manifold cannot have arbitrarily small
subgroups. By the work of Newman and Smith, it is sufficient to
prove the special case when the topological group is the $p$-adic
integers $\widehat{\mathbb{Z}}_p=\underset{\leftarrow}{lim} \{
\phi_n:\mathbb{Z}/p^{n+1} \mathbb{Z}\to \mathbb{Z}/p^n
\mathbb{Z}\}$, where $\phi_n$ is the $mod \ p^n$ mapping.

\bigskip
C.T.Yang showed that for any $n$-manifold $M^n$, the orbit space
$M/\widehat{\mathbb{Z}}_p$ must have integral cohomological
dimension $dim_{\mathbb{Z}} M/\widehat{\mathbb{Z}}_p=n+2$, and
more generally, for any locally compact Hausdorff cohomological
$n$-dimensional space X supporting an effective
$\widehat{\mathbb{Z}}_p$ action, $dim_{\mathbb{Z}}
X/\widehat{\mathbb{Z}}_p\leq n+3$. For locally compact finite
dimensional metric space $X$, $dim_{Z}\ X=dim\ X$, hence Yang's
results imply that if $X$ is metrizable, then $dim
X/\widehat{\mathbb{Z}}_p \leq n+3$ or $dim
X/\widehat{\mathbb{Z}}_p =\infty$.

\bigskip Although no effective action of $p$-adic group on
manifolds has yet been constructed, there do exist $p$-adic group
actions on Menger manifolds. Let $X$ be a $\mu^n$-manifold. Then
there are three constructions of $\widehat{\mathbb{Z}}{}_p$
actions on $X$.

\noindent (1) Every compact $0$-dimensional metrizable group G
acts effectively on $X$ so that $dim\ X/G=dim\ X$. Hence $p$-adic
group acts on Menger compacta. This is done by Dranishnikov \cite
{Dr2} and Mayer and Stark \cite {MS} using $Pasynkov's$ partial
product description of $\mu^n$. K. $Sakai$ \cite {Sa} has another
construction.

\noindent (2) $\widehat{\mathbb{Z}}{}_p$ acts freely on $X$ so
that $dim X/\widehat{\mathbb{Z}}{}_p=n+1$. This example depends on
the work of Bestvina, Edwards, Mayer and Stark \cite {MS}.

\noindent (3) $\widehat{\mathbb{Z}}{}_p$ acts on $X$ so that $dim
X/\widehat{\mathbb{Z}}{}_p=n+2$. This is done by Mayer and Stark
\cite {MS}, based on a construction of Raymond and Williams \cite
{RW} .

\bigskip In contrast, no example has been constructed where
$dim X/\widehat{\mathbb{Z}}_p= \infty$. A.N. Dranishnikov and
J.E.West proved that any compact metric space Y is the orbit space
of a $G$-cation on a metric compactum $X$ with $dim_{\mathbb{Z}_p}
X=1$, where $G=\prod_{i=1}^\infty ({\mathbb{Z}}/p)$. But their
paper has an error in their lemma 15.

The goal of this paper is to fix their error and generalize their
group $G$ to $\widehat{\mathbb{Z}}_p$.

\section{The construction}
\noindent {\bf Notations:}

\noindent (1) For a given set $S$ of positive integers, we define
$G_S=\prod_{i\in S} Z/{p^i}$.

\noindent (2) Given a sequence $K=\{k_i\}_{i=1}^{\infty}$, we
define $K^*=\{k_{2n-1}\}_{n=1}^{\infty}$.

\noindent (3) If $B=\{k_i\mid i\in \Lambda\}$, where $\Lambda$ is
a subset of $N$ (the set of all positive integers), we define
$K-B=\{k_i\}_{i\notin \Lambda}$.

\noindent (4) For the sets of integers in this paper, we will
allow them to contain repeating elements, and the union is always
a disjoint union. If you like, you can think them as sequences,
finite or infinite.

\bigskip Given a prime $p$, an infinite sequence $K$ of
positive integers, and any finite dimensional simplicial complex
$L$, we will define $\widehat{L}(K)$.

Let's work on a simplex first. Let $L=\Delta^m$, and $L'$ be the
first barycentric subdivision of $L$. Let
$\widehat{L^{(0)}}=L^{(0)}=\widetilde{L}_{(0)}$,
$\widehat{L^{(1)}}=L^{(1)}{}'=\widetilde{L}_{(1)}$, $K_0=K_1=K$,
$B_0=B_1=\emptyset$, $N_0=N_1=1$. We will inductively construct
the following diagram:

$$\xymatrix{& s_1^{-1}(L^{(0)}{}') \ar [ld]_{r_0}
\ar @{^{(}->}[d] & s_2^{-1}(L^{(1)}{}') \ar [ld]_{r_1} \ar
@{^{(}->}[d]
& \cdots & s_m^{-1}(L^{(m-1)}{}') \ar [ld]_{r_{m-1}} \ar @{^{(}->}[d]\\
\widetilde{L}_{(0)} \ar [d]^{q_0} & \widetilde{L}_{(1)} \ar
[d]^{q_1} & \widetilde{L}_{(2)} \ar [d]^{q_2} &
\cdots & \widetilde{L}_{(m)} \ar [d]^{q_m}\\
\widehat{L^{(0)}}=\widehat{L^{(0)}} \ar [d]^{p_0} \ar @{^{(}->}[r]
& \widehat{L^{(1)}}  \ar [d]^{p_1} \ar @{^{(}->}[r]&
\widehat{L^{(2)}}  \ar [d]^{p_2} \ar @{^{(}->}[r] & \cdots \ar
@{^{(}->}[r] & \widehat{L}=\widehat{L^{(m)}} \ar [d]^{p_m} \\
L^{(0)}{}' \ar @{^{(}->}[r] & L^{(1)}{}' \ar @{^{(}->}[r]&
L^{(2)}{}' \ar @{^{(}->}[r] & \cdots \ar @{^{(}->}[r] & L' }$$

Suppose $s_k:\widetilde{L}_{(k)}\overset{q_k}{\rightarrow}
\widehat{L^{(k)}}\overset{p_k}{\rightarrow} L^{(k)}{}'$ have been
constructed, such that:

\noindent (k.1) for each $(k+1)$-simplex $\Delta$ of $L^{(k+1)}$,
$s_k^{-1}(\Delta^{(k)})$ is a closed oriented $k$-manifold,

\noindent (k.2) the map $r_{k-1}:s_k^{-1}(L^{(k-1)}{}') \to
\widetilde{L}_{(k-1)}$ is a $N_{k-1}$-$1$ covering map,

\noindent (k.3) for any simplex $\sigma \in \Delta^{(k)}{}'$,
$s_k:s_k^{-1}(\sigma)\to \sigma$ is a $\prod_{i=1}^k N_i$ to $1$
covering map,

\noindent (k.4) for any $(k)$-simplex $\Delta$ of $L^{(k)}$,the
induced map $H^1(s_{k-1}^{-1}(\Delta^{(k-1)}),\mathbb{Z}/p) \to
H^1(s_{k}^{-1}(\Delta^{(k-1)}),\mathbb{Z}/p)$ is trivial,

\noindent (k.5)
$\pi_1(s_{k}^{-1}(\Delta^{(k-1)}))/\pi_1(s_{k-1}^{-1}(\Delta^{(k-1)}))\cong
G_{B_{K_{k-1}}}$,

\noindent (k.6) $G_{B_{k}}$ acts on $s_k^{-1}(L^{(k-1)}{}')$ with
$\widetilde{L}_{(k-1)}$ as orbit space. It is a free
 simplicial action,

\noindent (k.7) $\prod_{i=1}^k G_{B_k}$ acts on
$s_k^{-1}(\Delta^{(k)})$ simplicially, and the orbit space is
$\Delta^{(k)}$. Hence $\prod_{i=1}^k G_{B_k}$ acts on
$\widetilde{L}_{(k)}$ simplicially, and the orbit space is
$L^{(k)}{}'$,

\noindent (k.8) The map $q_k:\widetilde{L}_{(k)}-
s_k^{-1}(L^{(k-1)}{}')\longrightarrow
\widehat{L^{(k)}}-\widehat{L^{(k-1)}}$ is a homeomorphism.

\begin{lem}
$H^1(s_k^{-1}(\Delta^{(k)}))$ is a free abelian group. There
exists a regular covering space $E$ of $s_k^{-1}(\Delta^{(k)})$,
such that
$$\pi_1(E)/\pi_1(s_{k}^{-1}(\Delta^{(k)}))\cong
G_{B_{K_{k+1}}}$$ Hence $G_{B_{k+1}}$ acts on $E$ as the deck
transformation group.
\end{lem}
\begin{proof}
$s_{k}^{-1}(\Delta^{(k)})$ is an oriented $k$-dimensional closed
manifold, hence their $H_{k-1}$ is a free abelian group, so is its
$H^1$.

Let $G$ denote $\pi_1(s_{k}^{-1}(\Delta^{(k)}))$, $G^1$ be the
commutator subgroup of $G$, $\{g_1,g_2,\cdots,g_l\}$ be the set of
free generators of $H_1(s_k^{-1}(\Delta^{(k)}))$, and $T$ be the
torsion subgroup of $H_1$. Let $f$ denote the natural map
$\pi_1\to \pi_1/G^1\cong H_1$. Denote the first $l$ integers in
$K_k$ by $n_1,n_2,\cdots n_l$. Define $B_{k+1}=\{n_1,n_2,\cdots
n_l\}$, $N_{k+1}=N_k\times \prod_{i=1}^l p^{n_i}$. Let $G'$ denote
the subgroup of $\pi_1$ generated by $f^{-1}(\{g_1^{(n_{1})},
g_2^{(n_{2})}, \cdots, g_l^{(n_{l})}\}\cup T)$. ($G^1\subset G'$)
It is a normal subgroup of $G$. $G/G'\cong \prod_{i=1}^l
Z/(p^{n_i})=G_{B_{k+1}}$. There is a canonical regular covering
space of $s_k^{-1}(\Delta^{(k)})$, denoted as $E$, corresponding
to the group $G'$. $E$ is a closed oriented $k$-manifold.
$r_k:E\to s_k^{-1}(\Delta^{(k)})$ is a $|G:G'|=\prod_{i=1}^l
p^{n_i}$ covering. (This serves the proof for (k+1.2) and
(k+1.3)), and $G/G'\simeq G_{B_{k+1}}$ acts on $E$ as the deck
transformation group, with $s_k^{-1}(\Delta^{(k)})$ as orbit
space. (This serves the proof for (k+1.5) and (k+1.6)).
\end{proof}

\bigskip For any $(k+1)$-simplex $\Delta$ of $L^{(k+1)}$, we get
the covering space $E$. Let $cE$ be the cone of $E$. It is a
contractible oriented $(k+1)$-manifold with boundary. There is a
natural projection $s_{k+1}:cE\to \Delta$, which is a branch
covering. $G_{B_{k+1}}$ acts on $cE$ in the obvious way. (This
serves the proof for (k+1.7)). By (k.3), for each $\sigma\in
\Delta'$, $s_{k+1}: (s_{k+1})^{-1}(\sigma)\to \sigma$ is
$\prod_{i+1}^{k+1} N_i$ to $1$ covering map. (This proves
(k+1.3)). Let $K_{k+1}=K_k-B_{k+1}$.

\bigskip\noindent {\bf Case 1:} If $k+2\leq m$, then for any given a
$(k+2)$-simplex $\Delta$ of $L^{(k+2)}$, it has $k+3$ faces, say
$\Delta_1,\Delta_2, \cdots, \Delta_{k+3}$. For them, using the
above construction, we get cones $c_1E_1,c_2E_2,\cdots,
c_{k+3}E_{k+3}$.

For any $i\neq j$, $E_i$ is homeomorphic to $E_j$. $\Delta_i$
meets $\Delta_j$ along one $k$-simplex $\sigma$. We pick one
homeomorphism $h_{ij}:E_i\to E_j$, such that the following diagram
commute:
$$\xymatrix{E_i \ar [r]^{h_{ij}} \ar [d]^{r_k}&E_j \ar [d]^{r_k'}
 \\ s_k^{-1}(\Delta_i^{(k)}) \ar [r] \ar [d]^{s_k} &
s_k^{-1}(\Delta_j^{(k)}) \ar [d]^{s_k'}\\
\Delta_i^{(k)}  \ar [r]& \Delta_j^{(k)}}$$ and
$h_{ij}(r_ks_k)^{-1}(\sigma)=(r_k's_k')^{-1}(\sigma)$.

Denote the projection $E\to \Delta^{(k)}$ by $Pr$. Given
$\sigma_i\in E_i$ and $\sigma_j\in E_j$, if
$Pr(\sigma_1)=Pr(\sigma_2)=\sigma$ and
$h_{ij}(\sigma_1)=\sigma_2$, we glue the $(k+2)$-simplexes
$c_i\sigma_1,c_j\sigma_2$ of $c_iE_i,c_jE_j$ along $\sigma_1$ and
$\sigma_2$ via $h_{ij}$. Do the same for all the pairs $E_i$ and
$E_j$. It can be proved easily that each $\sigma$ is glued to
exactly one other $\sigma'$. So we get a $(k+1)$-manifold. (This
proves (k+1.3)).

Do the same for all $(k+2)$-simplexes of $L^{(k+2)}$, we get a
$k+1$ dimensional complex, denote it as $\widetilde{L}_{(k+1)}$.

\bigskip\noindent {\bf Case 2:} If $k+2>m$, then $k+1=m$.
For the unique $m$-simplex $\Delta$, we get a cone $cE$. It is the
$\widetilde{L}_{(m)}$.

\bigskip Now we want to have the following diagram:

$$\xymatrix{&s_{k+1}^{-1}(L^{(k)}{}') \ar [ld]
\ar @{^{(}->}[d]  \\ \widetilde{L}_{(k)} \ar [d] &
\widetilde{L}_{(k+1)} \ar [d] \\
\widehat{L^{(k)}} \ar [d] \ar @{^{(}->}[r]& \widehat{L^{(k+1)}}\ar [d]\\
L^{(k)}{}' \ar @{^{(}->}[r] & L^{(k+1)}{}'}$$

The only thing missing here is $\widehat{L^{(k+1)}}$. Since we
have the map $r_k\circ p_k:s_{k+1}^{-1}(L^{(k)}{}')
\overset{r_k}{\longrightarrow} \widetilde{L}_{(k)}
\overset{p_k}{\longrightarrow} \widehat{L^{(k)}}$, we can define
$\widehat{L^{(k+1)}}$ as the following: in
$\widetilde{L}_{(k+1)}$, take any $k$-simplexes
$\sigma_1,\sigma_2$, if $r_k\circ p_k(\sigma_1)=r_k\circ
p_k(\sigma_2)$, then we glue them together. (This proves (k+1.8),
and also serves the proof for (k+1.7)). Hence we get an inclusion:
$\widehat{L^{(k)}}\to \widehat{L^{(k+1)}}$. This finishes the
construction for the step $k+1$.

Let $K'=\cup_{i=1}^m B_i$. We have $G_{K'}$ acts on $\widehat{L}$
with $L'$ as orbit space.

\bigskip For general $m$-complex, we define
$\widehat{L}(K)$ to be the pull back:

$$\xymatrix{ \widehat{L} \ar[d]^p\ar[rr]^{homeo}
&&\widehat{L}{}'\ar[d] \ar[rr]&& \widehat{\Delta}\ar[d] \\ L
\ar[rr]^{1st}_{subdivsion}  && L' \ar[r] & {\triangle}^m \ar
[r]^{homeo}& ({\triangle}^m){}' }$$

\begin{lem}
For each $(k+1)$-simplex $\Delta$ of $L^{(k+1)}$, the induced map
$q_k^*:H^1(s_k^{-1}(\Delta^{(k)}),\mathbb{Z}/p) \to
H^1(E=s_{k+1}^{-1}(\Delta^{(k)}),\mathbb{Z}/p)$ is trivial.
\end{lem}
\begin{proof}
By the universal coefficient theorem for homology, we have: $0\to
H_1(C)\otimes \mathbb{Z}/p \to H_1(C;\mathbb{Z}/p) \to
Tor(H_0(C),\mathbb{Z}/p)\to 0$. In the case $H_0(C)$ is free, we
have $H_1(C)\otimes \mathbb{Z}/p\cong H_1(C,\mathbb{Z}/p)$.

Using the relation between $\pi_1$ and $H_1$, we get the following
commutative diagram:
$$\xymatrix{\pi_1(E) \ar [r] \ar [d] & H_1(E) \ar [d] \ar [r]^{\otimes
\mathbb{Z}/p}&H_1(E;\mathbb{Z}/p) \ar [d]
\\ \pi_1(s_k^{-1}(\Delta^{(k)})) \ar [r]  &
H_1(s_k^{-1}(\Delta^{(k)})) \ar [r]^{\otimes
\mathbb{Z}/p}&H_1(s_k^{-1}(\Delta^{(k)});\mathbb{Z}/p) }$$ The
image of $\pi_1(E)$ in $\pi_1(s_k^{-1}(\Delta^{(k)}))$ is $G'$,
while $G'$'s image in $H_1(s_k^{-1}(\Delta^{(k)});\mathbb{Z}/p)$
is $0$. The map $\pi_1(E)\to H_1(E) \to H_1(E;\mathbb{Z}/p)$ is
onto. Hence we know the map $q_k{}_*:H_1(E;\mathbb{Z}/p) \to
H_1(s_k^{-1}(\Delta^{(k)});\mathbb{Z}/p)$ is $0$.

For any chain complex $C$, by the universal coefficient theorem
for cohomology, we have: $0\to Ext(H_0(C),\mathbb{Z}/p)\to
H^1(C;\mathbb{Z}/p)\to Hom(H_1(C),\mathbb{Z}/p)\to 0$. In the case
$H_0(C)$ is free, we have $H^1(C;\mathbb{Z}/p)\cong
Hom(H_1(C),\mathbb{Z}/p)$. And we have the following commutative
diagram:
$$\xymatrix{H^1(E;\mathbb{Z}/p)\ar [r]  & Hom(H_1(E),\mathbb{Z}/p)
\\ H^1(s_k^{-1}(\Delta^{(k)});\mathbb{Z}/p) \ar [r] \ar [u]^{q_k^*} &
Hom(H_1(s_k^{-1}(\Delta^{(k)})),\mathbb{Z}/p)  \ar
[u]^{(q_k{}_*)^*}}$$

Since $q_k{}_*=0$, we know $q_k^*=0$.
\end{proof}

\begin{lem}
Given any space $M$ and any commutative diagram:
$$\xymatrix{& s_1^{-1}(L^{(0)}{}') \ar [ld]_{r_0}
\ar @{^{(}->}[d] & s_2^{-1}(L^{(1)}{}') \ar [ld]_{r_1} \ar
@{^{(}->}[d] & \cdots  & s_m^{-1}(L^{(m-1)}{}') \ar [ld]_{r_{m-1}}
\ar @{^{(}->}[d]\\ \widetilde{L}_{(0)} \ar [d]^{f_0} &
\widetilde{L}_{(1)} \ar [d]^{f_1} & \widetilde{L}_{(2)} \ar
[d]^{f_2}& \cdots &
\widetilde{L}_{(m)} \ar [d]^{f_m}\\
M \ar [r]^{id} & M \ar [r]^{id}& M \ar [r]^{id}& \cdots \ar
[r]^{id}& M}$$ it uniquely determines a map $f':\widehat{L}\to M$
such that, for any $x\in \widetilde{L}_{(k)}-
s_k^{-1}(L^{(k-1)}{}')$, $f'(q_k(x))=f_k(x)$, where $q_k$ is the
projection $\widetilde{L}_{(k)}\to \widehat{L^{k}}$.
\end{lem}
\begin{proof}
For any $x\in \widehat{L}$, there is a $k\geq 0$, such that $x\in
\widehat{L^{(k)}}-\widehat{L^{(k-1)}}$. The map
$q_k:\widetilde{L}_{(k)}- s_k^{-1}(L^{(k-1)}{}')\longrightarrow
\widehat{L^{(k)}}-\widehat{L^{(k-1)}}$ is a homeomorphism. We
define $f'(x)=f_k(q_k^{-1}(x))$. This well defines the map
$f':\widehat{L}\to M$.
\end{proof}
\begin{lem}
Given any subcomplex $K$ of $L$, integers $l>m\geq k\geq 0$, any
map $f:\widehat{K}\cup \widehat{L^{(k)}}\to B^l(\mathbb{Z}/p)$,
there exists an extension $\widetilde{f}:\widehat{L}\to
B^l(\mathbb{Z}/p)$.
\end{lem}
\begin{proof}
\noindent Without lose of generality, we suppose $L=\Delta^m$.
Since $B^l(\mathbb{Z}/p)$ is path connected, if in addition $k=0$,
the map $f$ can be extended to an $f:\widehat{K}\cup
\widehat{L^{(1)}}\to B^l(\mathbb{Z}/p)$. The following diagram
defines $\widetilde{f}:q_2^{-1}(\widehat{K^{(2)}})\cup
s_2^{-1}(L^{(1)}{}') \to B^l(\mathbb{Z}/p)$.

$$\xymatrix{&s_1^{-1}(L^{(0)}{}') \ar [ld]
\ar @{^{(}->}[d] & s_2^{-1}(L^{(1)}{}') \ar [ld] \ar @{^{(}->}[d]
& \cdots & s_m^{-1}(L^{(m-1)}{}') \ar [ld]_{r_{m-1}} \ar
@{^{(}->}[d]\\ \widetilde{L}_{(0)} \ar [d] & \widetilde{L}_{(1)}
\ar [d] & \widetilde{L}_{(2)} \ar [d] & \cdots &
\widetilde{L}_{(m)}
\ar [d] \\
\widehat{L^{(0)}}=L^{(0)} \ar @{^{(}->}[r] & \widehat{L^{(1)}} \ar
@{^{(}->}[r]& \widehat{L^{(2)}} \ar @{^{(}->}[r]& \cdots \ar
@{^{(}->}[r] & \widehat{L}=\widehat{L^{(m)}}}$$

We will construct the map by induction. Suppose we have
$\widetilde{f}:q_{k+1}^{-1}(\widehat{K^{(k+1)}})\cup
s_{k+1}^{-1}(L^{(k)}{}') \to B^l(\mathbb{Z}/p)$. Take any
$(k+1)$-simplex $\Delta\in L$, where $f$ is not defined on
$\widehat{\Delta }$, $\widetilde{f}$ is defined only on
$s_{k+1}^{-1}(\Delta^{(k)})$.

$$\xymatrix{cE \ar @{-->}[rd] \\
s_{k+1}^{-1}(\Delta^{(k)})=E \ar [d] \ar @{-->}[r]\ar
[rd]^{\overline{f}} \ar @{^{(}->}[u] &
E^l(\mathbb{Z}/p) \ar [d] \\
s_k^{-1}(\Delta^{(k)})\ar @{^{(}->}[r]& B^l(\mathbb{Z}/p)
 }$$

Since the map $r_k:E\to s_k^{-1}(\Delta^{(k)})$ induces a trivial
map between $H^1(\mathbb{Z}/p)$, the map $\widetilde{f}:E\to
(B^l(\mathbb{Z}/p))$ is null homotopic. Hence $\widetilde{f}$
lifts to a map $E \to E^l(\mathbb{Z}/p)$. On the other hand,
$\pi_i(E^l(\mathbb{Z}/p))=0$ for $i=1,2,\cdots, m$, we get an
extension $\widetilde{f}:cE\to E^l(\mathbb{Z}/p)\to
B^l(\mathbb{Z}/p)$. Do this for all the $(k+1)$-simplexes of $L$,
we have a $\widetilde{f}:\widetilde{L}_{(k+1)} \to
B^l(\mathbb{Z}/p)$.

By induction, we get a map $\widetilde{L}_{(m)} \to
B^l(\mathbb{Z}/p)$, and we have the commutative diagram:

$$\xymatrix{&\pi_1^{-1}(L^{(0)}{}') \ar [ld]
\ar @{^{(}->}[d] & \pi_2^{-1}(L^{(1)}{}') \ar [ld] \ar
@{^{(}->}[d] & \cdots   & s_m^{-1}(L^{(m-1)}{}') \ar [ld] \ar
@{^{(}->}[d]\\ \widetilde{L}_{(0)} \ar [d] & \widetilde{L}_{(1)}
\ar [d] & \widetilde{L}_{(2)} \ar [d]& \cdots &
\widetilde{L}_{(m)} \ar [d]\\
B^l(\mathbb{Z}/p) \ar [r]^{id} & B^l(\mathbb{Z}/p) \ar [r]^{id}&
B^l(\mathbb{Z}/p) \ar [r]^{id}& \cdots \ar [r]^{id} &
B^l(\mathbb{Z}/p)}$$

This defines the map $\widehat{L}\to B^l(\mathbb{Z}/p)$.
\end{proof}

\bigskip Given a finite dimensional simplicial complex $L$ and an
infinite sequence $K$ of positive integers, we can construct a
space $X$ as following:

\noindent Let $K_1=K^*$, $L_{1}=\widehat{L}(K_1)$, $B_1=K_1'$,
$K_2=(K-B_1')^*$. Inductively, let
$L_{k}=\widehat{L_{k-1}}(K_{k})$, $B_k=B_{k-1}\cup K_k'$,
$K_{k+1}=(K-B_k')^*$. (The way I define the $K_i$'s satisfies tow
conditions: (1) each $K_i$ is an infinite sequence, (2) $\sqcup
B_i= K$. )

\noindent Let $X$ be the inverse limit of $p_{k+1}:L_{k+1}\to
L_k$. The group $G=\prod_{i\in K} Z/{p^i}$ acts on $\overline{L}$,
and the quotient space is $L$.

\begin{thm}
For any compact metric space $L$, and any prime $p$, $L$ is the
orbit space of a $G$-action on a metric compactum $X$ with
$dim_{\mathbb{Z}_p} X = 1$, where $G$ is a compact group.
\end{thm}

\begin{proof}
If $L$ is a finite dimensional simplicial complex, the $X$ is the
space we construct above. .

To prove that $dim_{\mathbb{Z}/p}X\leq 1$, it is enough to show
that the restriction map $i^*: H^1(X;\mathbb{Z}_p) \rightarrow
H^1(A;\mathbb{Z}_p)$ is surjective for any closed subset $A\subset
X$. It is equivalent to show that any map $\phi : A \rightarrow
K(\mathbb{Z}_p;1)$ extends to $X$. We can use of the fact that
$K(\mathbb{Z}_p;1)$ is a neighborhood retract to relate $\phi$ to
maps of $|L_i|$ into $K(\mathbb{Z}_p;1)$. Consider the infinite
mapping cylinder $M_\infty$ of the inverse system $\{p_i:
|L_{i-1}|\leftarrow |L_i|\}_{i=3}^\infty$, that is, let $M_i$ be
the mapping cylinder of $p_i$, and for each $i$, identify the copy
of $|L_{i-1}|$ in $M_i$ with the domain end of $M_{i-1}$. Let
$C_j= \cup_{i=1}^j M_i$. Then the collapses generate an inverse
system $C_{i-1} \leftarrow C_j$ containing $|L_{j-1}|\leftarrow
|L_j|$ with inverse limit $Z=M_\infty \cup X$. Now $A\subset Z$
and $\phi:A \rightarrow K(\mathbb{Z}_p;1)$ extends to a map
$\Phi:U\rightarrow K(\mathbb{Z}_p;1)$ where $U$ is an open
neighborhood of $A$ in $Z$.

The collapses of the mapping cylinders $M_i$ are connected to the
identities of the $M_i$ by one parameter families of retractions
which extend to give a one-parameter family of retractions $r_t:
Z\rightarrow Z$ of $Z$, $0\leq t\leq 1$ with $r_1$ the identity
and $r_0$ the retraction onto $|L_1|$ and such that for
$t_n=(n-1)/n$, $r_{t_n}$ is the projection onto $C_n$. By
compactness, there is an $n$ such that $r_t(A)\subset U$ for all
$t\geq t_{n-1}$, so Borsuk's homotopy extension theorem implies
that $\phi$ extends to $X$ provided that $\Phi\circ r_n$ extends
to$|L_n|$. For $n$ sufficiently large, the union $T$ of all
subcomplexes of $L_n$ of the form $p^{-1}_n (\sigma)$ for
simplices $\sigma$ of $|L_{n-1}|$ that $\sigma\cap
r_{t_n}(A)\subset U$. Let $\Psi=\Phi|_T:T \rightarrow
K(\mathbb{Z}_p;1)$. By lemma, $\Psi$ can be extend
$\Psi:|L_{n}|\to K(\mathbb{Z}_p;1)$. Thus, $\phi$ extends to $X$,
hence $dim_{\mathbb{Z}_p}X \leq 1$.

To prove the general case, we first construct a space $X$ over the
Hilbert cube as follows:

For $n=1$, let $\widetilde{X}_1=I$, $K$ be any sequence of
positive integers. $K_1=K^*$, $B_1=\emptyset$. For each $n$
consider the projection $I^{n+1}\rightarrow I^n$ and the following
diagram:

$$\xymatrix{\widetilde{X}{}_n \ar[d]^{q_n} & \widehat{X}{}_{n+1}
\ar[l]_{p_{n+1}'} \ar[d]^{p_n'{}'} && \widetilde{X}{}_{n+1}
\ar[ll]_{p_{n+1}}\ar[dll]_{q_{n+1}}\\ I^n & I^{n+1}\ar[l] && }$$

where $q_n$ is assumed by induction, $\widehat{X}{}_{n+1}$ is a
pull back,
$\widetilde{X}_{n+1}=\widehat{\widehat{X}{}_{n+1}}(K_{n+1})$,
$K_{n+1}=(K-B_{n}')^*$, $B_{n+1}=B_n\cup K_{n+1}'$.  Let:
$q_{n+1}=p_{n+1}\circ p_{n}'' : \widetilde{X}_{n+1} \rightarrow
I^{n+1}$. Let $G_n=G_{K_n'}\times G_{n-1}$.

By induction, we have an action of $G_n$ on $\widetilde{X}{}_n$,
with orbit map $q_{n}$, so we get an induced action of $G_{n+1}$
on $\widetilde{X}_{n+1}$ with orbit map $q_{n+1}$.
$dim_{\mathbb{Z}_p}\widetilde{X}_{n+1}= 1$. Now form the inverse
limits:

$$\xymatrix{
\widetilde{X}{}_1 \ar[d]^{q_1} & \widetilde{X}{}_2 \ar[l]
\ar[d]^{q_2} & \widetilde{X}{}_3 \ar[l]\ar[d]_{q_3} & \cdots
\ar[l] & X \ar [l] \ar [d]^{q} \\I^1 & I^2\ar[l] &
I_3\ar[l]&\cdots \ar[l] & I^{\infty} \ar [l] }$$

We now have an action of infinite product $\prod G_{K_n'} \cong
G_K$ on $X$ with orbit map $q$.

Using a similar construction as before, i.e. the infinite mapping
cylinder, we can prove $dim_{\mathbb{Z}_p}X \leq 1$. Actually,
each $dim_{\mathbb{Z}_p}\widetilde{X}{}_k=1$, so the extension
clearly exists.

To complete the proof, let $L$ be any compact metric space, embed
$L$ in $I^\infty$, and let $X=q^{-1}(L)$.

\end{proof}
\bigskip\noindent {\bf Note:}

\noindent (1) If we choose $K=\{1,1,1,\cdots\}$, then $G=\prod
\mathbb{Z}/p$. This fixes the error of A.N.Dranishnikov and
J.E.West's paper.

\noindent (2) If we choose $K=\{1,2,3,\cdots\}$, then
$G=\prod_{i=1}^{\infty} \mathbb{Z}/p^i$. On the other hand,
$\widehat{\mathbb{Z}}{}_p$ is a subgroup of $G$, so we get a
$\widehat{\mathbb{Z}}{}_p$ group action on $\widehat{L}$. I guess
the quotient space has dimension infinity, but I haven't proven it
yet.

\begin{thm}
Given any integer $n\geq 2$, any prime $p$, any sequence $K$ of
positive integers, there is an action of $G=G_K$ on a compact,
two-dimensional, metric space $X$ such that $dim X/G=n$, moreover,
$dim X\times X=0$.
\end{thm}
\begin{proof}
This follows the same proof as the proof of Theorem A of [D].
\end{proof}
\begin{cor}
For any $p$, $K$ as above, there is an action of $G=G_K$ on a
compact, two-dimensional, metric space $X_{\infty}$ such that $dim
X_{\infty}/G=\infty$ and $dim X_{\infty}\times X_{\infty}=3$.
\end{cor}




\end{document}